\documentclass[11pt]{amsart}
\usepackage{amsmath,amssymb,amsthm,mathrsfs}
\usepackage[margin=1.1in]{geometry}
\usepackage{enumitem}
\usepackage[colorlinks=true,linkcolor=blue,citecolor=blue,urlcolor=blue]{hyperref}

\newtheorem{theorem}{Theorem}[section]
\newtheorem{proposition}[theorem]{Proposition}
\newtheorem{lemma}[theorem]{Lemma}
\newtheorem{corollary}[theorem]{Corollary}
\theoremstyle{definition}
\newtheorem{definition}[theorem]{Definition}
\newtheorem{example}[theorem]{Example}
\newtheorem{question}[theorem]{Question}
\theoremstyle{remark}
\newtheorem{remark}[theorem]{Remark}

\newcommand{\Z}{\mathbb{Z}}
\newcommand{\R}{\mathbb{R}}
\newcommand{\C}{\mathbb{C}}
\newcommand{\D}{\mathbb{D}}
\newcommand{\norm}[1]{\lVert #1\rVert}
\newcommand{\Int}{\operatorname{Int}}

\newcommand{\Ghat}{\widehat{G}}
\newcommand{\TE}{\textup{(T$_E$)}}

\makeatletter
\@namedef{subjclassname@2020}{\textup{2020} Mathematics Subject Classification}
\makeatother

\begin{document}

\title[Transitive and hypercyclic semigroups of composition operators]{Hypercyclic and topologically transitive\\ semigroups of composition operators \\ (revised and extended version)}

\author{Mohammad Javaheri}
\address{Department of Mathematics and Statistics, Siena University\\
515 Loudon Road, Loudonville, NY 12211, USA}
\email{mjavaheri@siena.edu}


\begin{abstract}
Let $X$ be a separable locally compact metric space, let $E\neq\{0\}$ be a normed space, and let $G$ be a semigroup of continuous self-maps of $X$. Every $\varphi\in G$ induces the composition operator $C_\varphi f=f\circ\varphi$ on the space $C(X,E)$ of continuous functions endowed with the compact-open topology. We characterize the topological transitivity of the semigroup $\Ghat=\{C_\varphi:\varphi\in G\}$ by a condition that combines the run-away property of $G$ with an approximate form of injectivity, and we show that, when $E$ is separable, $\Ghat$ is topologically transitive if and only if it is hypercyclic. In particular, if every element of $G$ is one-to-one and $G$ is run-away, then $\Ghat$ is topologically transitive, and hypercyclic when $E$ is separable. For a single map $\varphi$ we obtain: $C_\varphi$ is hypercyclic $\iff$ $C_\varphi$ is topologically transitive $\iff$ $\varphi$ is one-to-one and run-away, and $C_\varphi$ is mixing $\iff$ $\varphi$ is one-to-one and strongly run-away. If, moreover, $\varphi$ is an open map (which is automatic on manifolds by invariance of domain), then $C_\varphi$ is Devaney chaotic. On an open interval all of these properties are equivalent to $\varphi$ being one-to-one without fixed points.
\end{abstract}

\maketitle

\section{Introduction}

Throughout, $(X,d)$ is a metric space, $(E,\norm{\cdot})$ is a real or complex normed space with $E\neq\{0\}$, and $C(X,E)$ denotes the space of continuous functions $X\to E$ endowed with the compact-open topology. Since $E$ is metric, this is the topology of uniform convergence on compact sets: $f_n\to f$ if and only if $f_n|_K\to f|_K$ uniformly for every compact $K\subseteq X$. For $u\colon X\to E$ and $K\subseteq X$ we write $\norm{u}_K=\sup_{x\in K}\norm{u(x)}$, and $U(a,\delta)$ denotes the open ball of radius $\delta$ centred at $a$. Every continuous map $\varphi\colon X\to X$ induces a continuous linear operator
\[
C_\varphi\colon C(X,E)\to C(X,E),\qquad C_\varphi(f)=f\circ\varphi .
\]

Recall that a continuous map $T$ on a topological space $Y$ is called
\begin{enumerate}[label=(\roman*)]
\item \emph{topologically transitive} if for every pair of nonempty open sets $U,V\subseteq Y$ there exists $n\ge1$ such that $T^n(U)\cap V\neq\emptyset$;
\item \emph{mixing} if for every pair of nonempty open sets $U,V\subseteq Y$ one has $T^n(U)\cap V\neq\emptyset$ for all sufficiently large $n$;
\item \emph{hypercyclic} if there exists $y\in Y$ such that $\{T^n(y):n\ge1\}$ is dense in $Y$;
\item \emph{Devaney chaotic} if it is topologically transitive and its set of periodic points $\{y: T^ny=y \text{ for some } n\ge1\}$ is dense in $Y$.
\end{enumerate}
For a semigroup $\mathcal S$ of continuous self-maps of $Y$ (under composition), the action of $\mathcal S$ is \emph{topologically transitive} if for all nonempty open $U,V\subseteq Y$ there is $T\in\mathcal S$ with $T(U)\cap V\neq\emptyset$, and \emph{hypercyclic} if there exists $y\in Y$ whose orbit $\{T(y):T\in\mathcal S\}$ is dense in $Y$. For a single map $T$ these notions reduce to (i) and (iii) with $\mathcal S=\{T^n:n\ge1\}$.

If $G$ is a semigroup of continuous self-maps of $X$, then $\Ghat=\{C_\varphi:\varphi\in G\}$ is a semigroup of operators on $C(X,E)$, because $C_\varphi C_\psi=C_{\psi\circ\varphi}$; the $\Ghat$-orbit of $f$ is $\{f\circ\varphi:\varphi\in G\}$.

\begin{definition}\label{def:runaway}
A semigroup $G$ of self-maps of $X$ is \emph{run-away} if for every compact set $K\subseteq X$ there exists $\varphi\in G$ such that $\varphi(K)\cap K=\emptyset$. A map $\varphi\colon X\to X$ is \emph{run-away} if the semigroup $\{\varphi^n:n\ge1\}$ is run-away, and \emph{strongly run-away} if for every compact set $K\subseteq X$ there exists $N\ge1$ such that $\varphi^n(K)\cap K=\emptyset$ for all $n\ge N$.
\end{definition}

A run-away map has no periodic points (a finite orbit is compact and invariant), and no compact space admits a run-away map.

\subsection*{Background}
Dynamical properties of $C_\varphi$ on various subspaces of $C(X,E)$ have been studied extensively in connection with the topological properties of $\varphi$. For a holomorphic self-map $\varphi$ of the unit disc $\D$, $C_\varphi$ is hypercyclic on the space $H(\D)$ of holomorphic functions if and only if it is mixing, if and only if it is Devaney chaotic, and if and only if $\varphi$ is univalent with no fixed point in $\D$; see Shapiro~\cite{Shapiro} and the survey of B\`es~\cite{Bes}. More generally, for a simply connected domain $\Omega\subseteq\C$, Grosse-Erdmann and Mortini~\cite{GEM} proved that $C_\varphi$ is hypercyclic on $H(\Omega)$ if and only if $\varphi$ is univalent and run-away. For an open set $\Omega\subseteq\R^d$ and a real analytic self-map $\varphi$ of $\Omega$, Bonet and Doma\'nski~\cite{BD} proved that $C_\varphi$ is hypercyclic on the space $\mathscr A(\Omega)$ of real analytic functions if and only if $\varphi$ is injective, its Jacobian never vanishes on $\Omega$, and $\varphi$ is run-away. Since the inclusion $\mathscr A(I)\hookrightarrow C(I,\R)$ is continuous and has dense range (it contains the polynomials, which are dense by the Stone--Weierstrass theorem~\cite{Saxe}), the comparison principle~\cite[Prop.~1.19]{GEP} shows that hypercyclicity on $\mathscr A(I)$ implies hypercyclicity on $C(I,\R)$; Corollary~\ref{cor:interval} below gives a direct and complete answer on $C(I,E)$.

Continuous semigroups of composition operators have also been studied. Let $\rho\colon[0,\infty)\to(0,\infty)$ satisfy $\rho(\tau)\le Me^{\omega t}\rho(t+\tau)$ for all $\tau\ge0$, $t>0$ and some $M\ge1$, $\omega\in\R$, and let $C_{0,\rho}$ be the Banach space of continuous $u\colon[0,\infty)\to\C$ with $\rho(\tau)u(\tau)\to0$ as $\tau\to\infty$, normed by $\norm{u}=\sup_\tau|u(\tau)|\rho(\tau)$. Desch, Schappacher and Webb~\cite{DSW} proved that the translation semigroup $[T(t)u](\tau)=u(\tau+t)$ is hypercyclic on $C_{0,\rho}$ if and only if $\liminf_{t\to\infty}\rho(t)=0$, and Berm\'udez, Bonilla, Conejero and Peris~\cite{BBCP} proved that it is mixing if and only if $\lim_{t\to\infty}\rho(t)=0$; see also~\cite{BBCP,dLE,DSW} for weighted $L^p$ spaces and~\cite{CP} for complex sectors. The hypothesis on $\rho$ gives $\rho(t)\ge\rho(0)e^{-\omega t}/M$, so $\rho$ is bounded below on compact sets and the inclusion $C_{0,\rho}\hookrightarrow C([0,\infty),\C)$ is continuous with dense range; hence the translation semigroup is hypercyclic on $C([0,\infty),\C)$ as well. Example~\ref{ex:translations} recovers and generalizes this fact.

\subsection*{Main results}
The following condition turns out to characterize transitivity.

\begin{definition}\label{def:TE}
Let $G$ be a semigroup of continuous self-maps of $X$. We say that $G$ satisfies \TE{} if for every compact set $K\subseteq X$, every $g\in C(X,E)$ and every $\varepsilon>0$, there exists $\varphi\in G$ such that
\begin{enumerate}[label=(\alph*)]
\item $\varphi(K)\cap K=\emptyset$, and
\item $\norm{g(x)-g(y)}<\varepsilon$ whenever $x,y\in K$ and $\varphi(x)=\varphi(y)$.
\end{enumerate}
\end{definition}

Condition (b) is vacuous if $\varphi$ is injective on $K$. Thus \TE{} holds, for every $E$, as soon as
\begin{enumerate}[label=(I)]
\item for every compact $K\subseteq X$ there exists $\varphi\in G$ that is injective on $K$ and satisfies $\varphi(K)\cap K=\emptyset$;
\end{enumerate}
in particular, \TE{} holds if every element of $G$ is one-to-one and $G$ is run-away.

\begin{theorem}\label{thm:main}
Let $X$ be a separable locally compact metric space, $E\neq\{0\}$ a normed space, and $G$ a semigroup of continuous self-maps of $X$. Then:
\begin{enumerate}[label=\textup{(\arabic*)}]
\item The action of $\Ghat$ on $C(X,E)$ is topologically transitive if and only if $G$ satisfies \TE.
\item If $E$ is separable, the action of $\Ghat$ on $C(X,E)$ is hypercyclic if and only if it is topologically transitive, if and only if $G$ satisfies \TE.
\item If $E$ is not separable, the action of $\Ghat$ is never hypercyclic.
\end{enumerate}
Moreover, every hypercyclic function $f$ for $\Ghat$ has dense range $f(X)$ in $E$.
\end{theorem}

Part~(1) holds for an arbitrary metric space $X$ (see Proposition~\ref{prop:transitive}). Birkhoff's transitivity theorem~\cite{BM,GEP} gives the equivalence of transitivity and hypercyclicity for a single operator on a separable completely metrizable space without isolated points, but no such general principle is available for semigroup actions, and we do not assume $E$ to be complete. The following corollary contains the main theorem of the original version of this paper (where $E$ was assumed finite-dimensional).

\begin{corollary}\label{cor:injective}
Let $X$ be a separable locally compact metric space, $E\neq\{0\}$ a normed space, and $G$ a semigroup of continuous self-maps of $X$ satisfying \textup{(I)}; for example, every element of $G$ is one-to-one and $G$ is run-away. Then the action of $\Ghat$ on $C(X,E)$ is topologically transitive, and it is hypercyclic if $E$ is separable.
\end{corollary}

Injectivity of all elements of $G$ is not necessary; see Example~\ref{ex:subbiah}. For a single map we obtain a complete description.

\begin{theorem}\label{thm:single}
Let $X$ be a separable locally compact metric space, $E\neq\{0\}$ a normed space, and $\varphi\colon X\to X$ continuous.
\begin{enumerate}[label=\textup{(\arabic*)}]
\item $C_\varphi$ is topologically transitive on $C(X,E)$ if and only if $\varphi$ is one-to-one and run-away.
\item If $E$ is separable, $C_\varphi$ is hypercyclic if and only if it is topologically transitive.
\item $C_\varphi$ is mixing if and only if $\varphi$ is one-to-one and strongly run-away.
\end{enumerate}
\end{theorem}

\begin{theorem}\label{thm:chaos}
Let $X$ be a separable locally compact metric space, $E\neq\{0\}$ a normed space, and $\varphi\colon X\to X$ a continuous, one-to-one, open map that is strongly run-away. Then $C_\varphi$ is mixing and Devaney chaotic on $C(X,E)$, and hypercyclic if $E$ is separable.
\end{theorem}

By invariance of domain, openness is automatic on manifolds, which yields the following.

\begin{corollary}\label{cor:manifold}
Let $X$ be a second countable topological manifold without boundary all of whose components have the same dimension, and let $E\neq\{0\}$ be a normed space. For a continuous map $\varphi\colon X\to X$ the following are equivalent:
\begin{enumerate}[label=\textup{(\alph*)}]
\item $C_\varphi$ is mixing on $C(X,E)$;
\item $C_\varphi$ is mixing and Devaney chaotic on $C(X,E)$;
\item $\varphi$ is one-to-one and strongly run-away.
\end{enumerate}
\end{corollary}

\begin{corollary}\label{cor:interval}
Let $I\subseteq\R$ be an open interval, $E\neq\{0\}$ a normed space, and $\varphi\colon I\to I$ continuous. The following are equivalent:
\begin{enumerate}[label=\textup{(\alph*)}]
\item $C_\varphi$ is topologically transitive on $C(I,E)$;
\item $C_\varphi$ is mixing on $C(I,E)$;
\item $C_\varphi$ is Devaney chaotic on $C(I,E)$;
\item $\varphi$ is one-to-one and has no fixed points.
\end{enumerate}
If $E$ is separable, these are also equivalent to: $C_\varphi$ is hypercyclic on $C(I,E)$.
\end{corollary}

\begin{corollary}\label{cor:semiconj}
Let $X$ be a separable locally compact metric space, $E\neq\{0\}$ a normed space, and $\varphi\colon X\to X$ a continuous, one-to-one, open map. Suppose that there exist an open interval $J\subseteq\R$, a continuous map $\theta\colon X\to J$, and a continuous map $\psi\colon J\to J$ without fixed points such that $\theta\circ\varphi=\psi\circ\theta$. Then $C_\varphi$ is mixing and Devaney chaotic on $C(X,E)$, and hypercyclic if $E$ is separable.
\end{corollary}

The paper is organized as follows. Section~\ref{sec:prelim} collects preliminaries. Theorem~\ref{thm:main} and Corollary~\ref{cor:injective} are proved in Section~\ref{sec:main}, together with examples. Section~\ref{sec:single} treats single maps and Section~\ref{sec:chaos} treats chaos. Section~\ref{sec:remarks} contains further remarks and open questions.

\section{Preliminaries}\label{sec:prelim}

\begin{lemma}\label{lem:basis}
Let $X$ be a topological space and $E$ a normed space. For $f\in C(X,E)$, compact $K\subseteq X$ and $\varepsilon>0$, put
\[
B(f;K,\varepsilon)=\{h\in C(X,E):\norm{h-f}_K<\varepsilon\}.
\]
These sets form a basis of the compact-open topology of $C(X,E)$. Consequently, if $W\subseteq C(X,E)$ is open and $f\in W$, then there exist a compact set $K\subseteq X$ and $\varepsilon>0$ such that every $h\in C(X,E)$ with $\norm{h-f}_K<\varepsilon$ lies in $W$; in particular, every $h$ with $h|_K=f|_K$ lies in $W$.
\end{lemma}

\begin{proof}
For metric targets, the compact-open topology coincides with the topology of compact convergence, whose basic neighbourhoods are the sets $B(f;K,\varepsilon)$; see~\cite[Thm.~46.8]{Munkres}.
\end{proof}

\begin{remark}
The compact set $K$ in Lemma~\ref{lem:basis} depends on $f$; in general there is no compact $K$ that works for all $f\in W$ simultaneously. (For $X=\R$ and $W=\bigcup_{x\in\R}\{f:|f(x)|>1\}$, any compact $K$ fails.) In the proof of transitivity we therefore fix $f\in V$ and $g\in W$ \emph{before} choosing the compact sets.
\end{remark}

\begin{lemma}\label{lem:exhaustion}
Let $X$ be a separable locally compact metric space. Then:
\begin{enumerate}[label=\textup{(\alph*)}]
\item There are compact sets $L_1\subseteq L_2\subseteq\cdots$ with $L_j\subseteq\Int L_{j+1}$ for all $j$ and $\bigcup_j L_j=X$.
\item Every compact $K\subseteq X$ is contained in some $L_j$.
\item For all $j,l\ge1$ there is a finite set $J_{j,l}=\{a_1,\dots,a_k\}\subseteq L_j$, listed in a fixed order, with $L_j\subseteq\bigcup_{i=1}^kU(a_i,1/l)$.
\item If $Y$ is a topological space and $f\colon X\to Y$ is a function such that $f|_{L_j}$ is continuous for every $j$, then $f$ is continuous.
\end{enumerate}
\end{lemma}

\begin{proof}
(a) holds for every second countable locally compact Hausdorff space~\cite[Prop.~4.76]{Lee}. (b) The sets $\Int L_j$ form an increasing open cover of $K$, so $K\subseteq\Int L_j\subseteq L_j$ for some $j$. (c) follows from compactness of $L_j$. (d) If $O\subseteq Y$ is open, then $f^{-1}(O)\cap\Int L_j=(f|_{L_j})^{-1}(O)\cap\Int L_j$ is open in $X$, and $f^{-1}(O)=\bigcup_j\big(f^{-1}(O)\cap\Int L_j\big)$.
\end{proof}

By Lemma~\ref{lem:exhaustion}, $X$ is hemicompact and
\[
\rho(f,g)=\sum_{j=1}^\infty 2^{-j}\min\{1,\norm{f-g}_{L_j}\}
\]
is a metric on $C(X,E)$ inducing the compact-open topology. Thus $C(X,E)$ is a metrizable locally convex space (a Fr\'echet space if $E$ is a Banach space), and it has no isolated points.

We shall repeatedly use Dugundji's extension theorem~\cite{Dugundji}.

\begin{lemma}[Dugundji]\label{lem:dugundji}
Let $Y$ be a metric space, $A\subseteq Y$ closed, and $E$ a locally convex space (e.g.\ a normed space). Every continuous map $u\colon A\to E$ has a continuous extension $\tilde u\colon Y\to E$ with $\tilde u(Y)$ contained in the convex hull of $u(A)$.
\end{lemma}

For finite-dimensional $E$ this is just the Tietze extension theorem applied to coordinates. We also need explicit partitions of unity.

\begin{lemma}\label{lem:pou}
Let $Y$ be a metric space and let $W_1,\dots,W_r\subseteq Y$ be open sets with union $W$. With the convention $d(y,\emptyset)=1$, the functions
\[
\rho_i(y)=\frac{d(y,Y\setminus W_i)}{\sum_{j=1}^rd(y,Y\setminus W_j)},\qquad y\in W,
\]
are continuous on $W$, satisfy $\rho_i\ge0$ and $\sum_i\rho_i=1$ on $W$, and $\rho_i(y)>0$ if and only if $y\in W_i$.
\end{lemma}

\begin{proof}
For $y\in W$ some $d(y,Y\setminus W_j)$ is positive, since $Y\setminus W_j$ is closed. The rest is immediate.
\end{proof}

The next lemma is the key technical step: it constructs an approximate ``inverse composition'' for maps that are not injective.

\begin{lemma}\label{lem:approxinverse}
Let $X$ be a metric space, $E$ a normed space, $K\subseteq X$ compact, $\varphi\colon X\to X$ continuous, $g\in C(K,E)$ and $\varepsilon>0$. Suppose that $\norm{g(x)-g(y)}\le\varepsilon$ whenever $x,y\in K$ and $\varphi(x)=\varphi(y)$. Then there exists $F\in C(\varphi(K),E)$ such that
\[
\norm{F(\varphi(x))-g(x)}<2\varepsilon\qquad\text{for all }x\in K.
\]
If $\varphi$ is injective on $K$, one may take $F=g\circ(\varphi|_K)^{-1}$, so that $F\circ\varphi=g$ on $K$.
\end{lemma}

\begin{proof}
If $\varphi|_K$ is injective, it is a continuous bijection from a compact space onto the Hausdorff space $\varphi(K)$, hence a homeomorphism, and $F=g\circ(\varphi|_K)^{-1}$ works. In general, by uniform continuity of $g$ on $K$ choose $\eta>0$ with $\norm{g(x)-g(x')}<\varepsilon$ whenever $x,x'\in K$, $d(x,x')<\eta$. For $z\in\varphi(K)$ let $\Phi_z=K\cap\varphi^{-1}(z)$, a nonempty compact set.

\emph{Claim.} Each $z\in\varphi(K)$ has an open neighbourhood $O_z$ in $\varphi(K)$ such that every $x\in K$ with $\varphi(x)\in O_z$ satisfies $d(x,\Phi_z)<\eta$. Otherwise there are $x_n\in K$ with $\varphi(x_n)\to z$ and $d(x_n,\Phi_z)\ge\eta$; passing to a subsequence, $x_n\to x\in K$, and then $\varphi(x)=z$, i.e.\ $x\in\Phi_z$, while $d(x,\Phi_z)\ge\eta$ by continuity of $d(\cdot,\Phi_z)$, a contradiction.

By compactness, $\varphi(K)=O_{z_1}\cup\dots\cup O_{z_r}$. Choose $x_i\in\Phi_{z_i}$, let $\rho_1,\dots,\rho_r$ be the partition of unity of Lemma~\ref{lem:pou} on $Y=\varphi(K)$ subordinate to $\{O_{z_i}\}$, and set $F=\sum_{i}\rho_i\,g(x_i)$. Let $x\in K$ and $i$ with $\rho_i(\varphi(x))>0$. Then $\varphi(x)\in O_{z_i}$, so there is $x'\in\Phi_{z_i}$ with $d(x,x')<\eta$, whence $\norm{g(x)-g(x')}<\varepsilon$; and $x',x_i\in\Phi_{z_i}$ give $\norm{g(x')-g(x_i)}\le\varepsilon$. Hence $\norm{g(x)-g(x_i)}<2\varepsilon$ for every such $i$, and since $F(\varphi(x))$ is a convex combination of these $g(x_i)$, $\norm{F(\varphi(x))-g(x)}<2\varepsilon$.
\end{proof}

\section{Transitivity and hypercyclicity of semigroups}\label{sec:main}

\begin{proposition}[Necessity]\label{prop:necessity}
Let $X$ be a metric space, $E\neq\{0\}$ a normed space and $G$ a semigroup of continuous self-maps of $X$. If the action of $\Ghat$ on $C(X,E)$ is topologically transitive or hypercyclic, then $G$ satisfies \TE.
\end{proposition}

\begin{proof}
Fix a compact $K\subseteq X$, $g\in C(X,E)$ and $\varepsilon>0$.

\emph{Transitive case.} Choose $c\in E$ with $\norm{c}>\norm{g}_K+\varepsilon$ and consider the nonempty open sets $V=B(0;K,\varepsilon/2)$ and $W=B(g+c;K,\varepsilon/2)$, where $c$ also denotes the constant function. By transitivity there are $\varphi\in G$ and $f\in V$ with $f\circ\varphi\in W$. If $x\in K$ and $\varphi(x)\in K$, then $\norm{f(\varphi(x))}<\varepsilon/2$, whereas $\norm{f(\varphi(x))-g(x)-c}<\varepsilon/2$ forces $\norm{f(\varphi(x))}>\norm{c}-\norm{g}_K-\varepsilon/2>\varepsilon/2$; hence $\varphi(K)\cap K=\emptyset$. If $x,y\in K$ and $\varphi(x)=\varphi(y)$, then
\[
\norm{g(x)-g(y)}\le\norm{g(x)+c-f(\varphi(x))}+\norm{f(\varphi(y))-g(y)-c}<\varepsilon .
\]

\emph{Hypercyclic case.} Let $f$ be a hypercyclic function and choose $c\in E$ with $\norm{c}>\norm{f}_K+\norm{g}_K+\varepsilon$. Pick $\varphi\in G$ with $\norm{f\circ\varphi-g-c}_K<\varepsilon/2$. If $x\in K$ and $\varphi(x)\in K$, then $\norm{f(\varphi(x))}\le\norm f_K$, but also $\norm{f(\varphi(x))}>\norm c-\norm g_K-\varepsilon/2>\norm f_K$, a contradiction. Condition (b) follows exactly as before.
\end{proof}

\begin{proposition}[Sufficiency for transitivity]\label{prop:transitive}
Let $X$ be a metric space, $E$ a normed space and $G$ a semigroup of continuous self-maps of $X$ satisfying \TE. Then the action of $\Ghat$ on $C(X,E)$ is topologically transitive.
\end{proposition}

\begin{proof}
Fix $f_0\in V$ and $g_0\in W$. By Lemma~\ref{lem:basis} there are compact sets $K_1,K_2$ and $\varepsilon>0$ with $B(f_0;K_1,\varepsilon)\subseteq V$ and $B(g_0;K_2,\varepsilon)\subseteq W$. Put $K=K_1\cup K_2$ and let $\varphi\in G$ satisfy (a) and (b) of \TE{} for $K$, $g_0$ and $\varepsilon/2$. By Lemma~\ref{lem:approxinverse} (applied with $\varepsilon/2$) there is $F\in C(\varphi(K),E)$ with $\norm{F(\varphi(x))-g_0(x)}<\varepsilon$ for $x\in K$. Since $K$ and $\varphi(K)$ are disjoint compact sets, the function $u$ equal to $f_0$ on $K$ and to $F$ on $\varphi(K)$ is continuous on the closed set $K\cup\varphi(K)$; by Lemma~\ref{lem:dugundji} it extends to some $f\in C(X,E)$. Then $f=f_0$ on $K_1$, so $f\in V$, and for $x\in K_2$ we have $\norm{f(\varphi(x))-g_0(x)}=\norm{F(\varphi(x))-g_0(x)}<\varepsilon$; as $K_2$ is compact, $\norm{f\circ\varphi-g_0}_{K_2}<\varepsilon$ and $C_\varphi f\in W$.
\end{proof}

\begin{proposition}[Sufficiency for hypercyclicity]\label{prop:hypercyclic}
Let $X$ be a separable locally compact metric space, $E$ a separable normed space and $G$ a semigroup of continuous self-maps of $X$ satisfying \TE. Then the action of $\Ghat$ on $C(X,E)$ is hypercyclic.
\end{proposition}

\begin{proof}
Let $(L_j)$ and $J_{j,l}$ be as in Lemma~\ref{lem:exhaustion}, and let $\mathcal D\subseteq E$ be a countable dense set. Consider all triples $(j,l,b)$ with $j,l\ge1$ and $b=(b_1,\dots,b_k)\in\mathcal D^k$, where $k$ is the cardinality of $J_{j,l}$ (repetitions among the $b_i$ are allowed). There are countably many such triples; enumerate them as $(j_n,l_n,b^{(n)})$, $n\ge1$. For each $n$ write
\[
M_n=L_{j_n},\qquad J_{j_n,l_n}=(a_1,\dots,a_k),\qquad b^{(n)}=(b_1,\dots,b_k),
\]
let $\theta_1,\dots,\theta_k$ be the partition of unity of Lemma~\ref{lem:pou} subordinate to the open cover $\{U(a_i,1/l_n)\}_{i=1}^k$ of $M_n$, and define
\begin{equation}\label{eq:hn}
h_n(x)=\sum_{i=1}^k\theta_i(x)\,b_i,\qquad x\in M_n .
\end{equation}
By Lemma~\ref{lem:dugundji} we regard $h_n$ as an element of $C(X,E)$.

We construct inductively compact sets $K_0=\emptyset\subseteq K_1\subseteq K_2\subseteq\cdots$, functions $f_n\in C(K_n,E)$ and maps $\varphi_n\in G$ ($n\ge1$) such that
\begin{enumerate}[label=(\roman*)]
\item $L_n\subseteq K_n$ and $K_{n-1}\subseteq K_n$;
\item $f_n|_{K_{n-1}}=f_{n-1}$;
\item $\varphi_n(M_n)\subseteq K_n$ and $\norm{f_n(\varphi_n(x))-h_n(x)}<2/l_n$ for all $x\in M_n$.
\end{enumerate}
Suppose $K_{n-1}$ and $f_{n-1}$ have been constructed. Apply \TE{} to the compact set $C=K_{n-1}\cup M_n$, the function $h_n$ and $1/l_n$: there is $\varphi_n\in G$ with $\varphi_n(C)\cap C=\emptyset$ and $\norm{h_n(x)-h_n(y)}<1/l_n$ whenever $x,y\in M_n$ and $\varphi_n(x)=\varphi_n(y)$. Lemma~\ref{lem:approxinverse} gives $F_n\in C(\varphi_n(M_n),E)$ with $\norm{F_n\circ\varphi_n-h_n}_{M_n}<2/l_n$. Put
\[
K_n=K_{n-1}\cup L_n\cup\varphi_n(M_n).
\]
The compact sets $K_{n-1}$ and $\varphi_n(M_n)\subseteq\varphi_n(C)$ are disjoint, so the function equal to $f_{n-1}$ on $K_{n-1}$ and to $F_n$ on $\varphi_n(M_n)$ is continuous on their (closed) union; let $f_n\in C(K_n,E)$ be an extension of it given by Lemma~\ref{lem:dugundji}. Properties (i)--(iii) hold.

Define $f\colon X\to E$ by $f(x)=f_n(x)$ for $x\in K_n$. By (i) and (ii) this is well defined on $\bigcup_nK_n\supseteq\bigcup_nL_n=X$, and $f|_{L_n}=f_n|_{L_n}$ is continuous, so $f\in C(X,E)$ by Lemma~\ref{lem:exhaustion}(d). Moreover, by (iii), $\norm{f\circ\varphi_n-h_n}_{M_n}<2/l_n$ for all $n$.

Let $g\in C(X,E)$, $K\subseteq X$ compact and $\varepsilon>0$. By Lemma~\ref{lem:exhaustion}(b), $K\subseteq L_t$ for some $t$. Choose $m\ge1$ such that $2/m<\varepsilon/3$ and $\norm{g(x)-g(y)}<\varepsilon/3$ whenever $x,y\in L_t$ and $d(x,y)<1/m$. Write $J_{t,m}=(a_1,\dots,a_s)$ and choose $b_i\in\mathcal D$ with $\norm{b_i-g(a_i)}<\varepsilon/3$. There is $n$ with $(j_n,l_n,b^{(n)})=(t,m,(b_1,\dots,b_s))$; then $M_n=L_t$. For $x\in L_t$ let $I_x=\{i:d(x,a_i)<1/m\}$; then $\theta_i(x)=0$ for $i\notin I_x$, and for $i\in I_x$,
\[
\norm{b_i-g(x)}\le\norm{b_i-g(a_i)}+\norm{g(a_i)-g(x)}<2\varepsilon/3 .
\]
Hence $\norm{h_n(x)-g(x)}=\big\|\sum_{i\in I_x}\theta_i(x)(b_i-g(x))\big\|<2\varepsilon/3$ for $x\in L_t$, and
\[
\norm{f\circ\varphi_n-g}_K\le\norm{f\circ\varphi_n-h_n}_{L_t}+\norm{h_n-g}_{L_t}<\varepsilon/3+2\varepsilon/3=\varepsilon
\]
(the second supremum is attained on the compact set $L_t$). Thus the $\Ghat$-orbit of $f$ is dense.
\end{proof}

\begin{proof}[Proof of Theorem~\ref{thm:main}]
(1) follows from Propositions~\ref{prop:necessity} and~\ref{prop:transitive}, and (2) from Propositions~\ref{prop:necessity}, \ref{prop:transitive} and~\ref{prop:hypercyclic}. For (3) and the final assertion, let $f$ be hypercyclic, fix $x_0\in X$ and $e\in E$. Some $f\circ\varphi$ lies in $B(e;\{x_0\},\delta)$ for any $\delta>0$, so $\norm{f(\varphi(x_0))-e}<\delta$; hence $f(X)$ is dense in $E$. Since $X$ is separable, so is $f(X)$, and therefore $E$ is separable.
\end{proof}

\begin{proof}[Proof of Corollary~\ref{cor:injective}]
Condition (I) implies \TE{}, so Theorem~\ref{thm:main} applies.
\end{proof}

\begin{example}[Translation semigroups]\label{ex:translations}
Let $\Gamma\subseteq\R^d$ be an additive subsemigroup and let $X\subseteq\R^d$ be closed with $X+\Gamma\subseteq X$ (e.g.\ $X=\R^d$, or $X=[0,\infty)$ and $\Gamma\subseteq[0,\infty)$). Let $G=\{\tau_v:v\in\Gamma\}$, $\tau_v(x)=x+v$. All $\tau_v$ are injective, and $G$ is run-away if and only if $\Gamma$ is unbounded: if $K$ is compact and $|v|>\operatorname{diam}K$ then $(K+v)\cap K=\emptyset$, whereas if $|v|\le R$ for all $v\in\Gamma$, then $(B+v)\cap B\ne\emptyset$ for every $v\in\Gamma$, where $B=X\cap\overline{U}(x_0,R)$ for any fixed $x_0\in X$ (it contains $x_0$ and $x_0+v$). Hence, for separable $E\neq\{0\}$, the translation semigroup $\{C_{\tau_v}:v\in\Gamma\}$ is hypercyclic on $C(X,E)$ if and only if $\Gamma$ is unbounded. With $X=\Gamma=[0,\infty)$ this is the translation $C_0$-semigroup on the Fr\'echet space $C([0,\infty),E)$ discussed in the introduction.
\end{example}

\begin{example}[Non-injective generators]\label{ex:subbiah}
Let $I=(0,1)$. Since $I$ is homeomorphic to $\R$, a theorem of Subbiah~\cite{Subbiah} (see also Knichal~\cite{Knichal}) provides continuous maps $\eta,\mu\colon I\to I$ such that the semigroup $\langle\eta,\mu\rangle$ they generate is dense in the space $C(I,I)$ of continuous self-maps of $I$ with the compact-open topology. Let $f\colon I\to\R$ be a homeomorphism. Given $g\in C(I,\R)$, put $\vartheta=f^{-1}\circ g\in C(I,I)$ and choose $\lambda_i\in\langle\eta,\mu\rangle$ with $\lambda_i\to\vartheta$ compactly. If $K\subseteq I$ is compact, $\vartheta(K)$ is a compact subset of $I$, so for some $\delta>0$ the closed $\delta$-neighbourhood $K'$ of $\vartheta(K)$ is a compact subset of $I$, and $\lambda_i(K)\subseteq K'$ for $i$ large. Since $f$ is uniformly continuous on $K'$, $f\circ\lambda_i\to f\circ\vartheta=g$ uniformly on $K$. Hence the orbit of $f$ under the semigroup $\langle C_\eta,C_\mu\rangle$ is dense in $C(I,\R)$. A semigroup dense in $C(I,I)$ must contain non-injective maps, so injectivity of every element of $G$ is not necessary in Corollary~\ref{cor:injective}; by Theorem~\ref{thm:main}, $\langle\eta,\mu\rangle$ nevertheless satisfies \TE{} for $E=\R$.
\end{example}

\section{Single composition operators}\label{sec:single}

\begin{proof}[Proof of Theorem~\ref{thm:single}]
Let $G=\{\varphi^n:n\ge1\}$.

(1) If $\varphi$ is injective and run-away, then $G$ satisfies (I), hence \TE{}, and $C_\varphi$ is transitive by Proposition~\ref{prop:transitive}. Conversely, if $C_\varphi$ is transitive, then $G$ satisfies \TE{} by Proposition~\ref{prop:necessity}, so $\varphi$ is run-away. Suppose $\varphi(a)=\varphi(b)$ with $a\neq b$. Let $K=\{a,b\}$, fix $e\in E$ with $\norm e=1$, and let $g(x)=\frac{d(x,a)}{d(x,a)+d(x,b)}\,e$, so that $g(a)=0$, $g(b)=e$. Every $\varphi^n$ identifies $a$ and $b$, so (b) of \TE{} fails for $\varepsilon=1/2$, a contradiction. Thus $\varphi$ is injective.

(2) This is Theorem~\ref{thm:main}(2) for $G$.

(3) Suppose $\varphi$ is injective and strongly run-away, and let $V,W$ be nonempty open sets. With $K$ as in the proof of Proposition~\ref{prop:transitive}, choose $N$ with $\varphi^n(K)\cap K=\emptyset$ for $n\ge N$. Each $\varphi^n$ is injective, so it satisfies (a) and (b) of \TE{} for $K$ and any $g,\varepsilon$; the proof of Proposition~\ref{prop:transitive} gives $C_\varphi^n(V)\cap W=C_{\varphi^n}(V)\cap W\neq\emptyset$ for all $n\ge N$. Conversely, let $C_\varphi$ be mixing. Then it is transitive, so $\varphi$ is injective by (1). Given a compact $K$, take $c\in E$ with $\norm c>1$, $V=B(0;K,1/2)$ and $W=B(c;K,1/2)$. For $n\ge N$ there is $f\in V$ with $f\circ\varphi^n\in W$, and the argument of Proposition~\ref{prop:necessity} shows $\varphi^n(K)\cap K=\emptyset$ for all $n\ge N$.
\end{proof}

\begin{remark}
For a single operator and a Banach space $E$, the equivalence in Theorem~\ref{thm:single}(2) also follows from Birkhoff's theorem, since $C(X,E)$ is then a separable Fr\'echet space. Our proof is constructive and needs neither completeness of $E$ nor Baire category.
\end{remark}

\section{Chaotic composition operators}\label{sec:chaos}

\begin{proof}[Proof of Theorem~\ref{thm:chaos}]
By Theorem~\ref{thm:single}(3), $C_\varphi$ is mixing, hence transitive, and it is hypercyclic if $E$ is separable by Theorem~\ref{thm:single}(2). It remains to show that periodic points are dense. We may assume $X\neq\emptyset$; then $X$ is not compact, since it admits a run-away map.

Let $g\in C(X,E)$ and $K\subseteq X$ compact; by Lemma~\ref{lem:basis} it suffices to find $F\in C(X,E)$ with $F=g$ on $K$ and $F\circ\varphi^p=F$ for some $p\ge1$. By Lemma~\ref{lem:exhaustion} there is a compact $K'$ with $K\subseteq\Int K'$. Let
\[
\chi(x)=\frac{d(x,X\setminus\Int K')}{d(x,X\setminus\Int K')+d(x,K)},\qquad h=\chi g .
\]
(the denominator never vanishes, and $X\setminus\Int K'\neq\emptyset$ because $X$ is not compact). Then $h\in C(X,E)$, $h=g$ on $K$, and $h=0$ outside $K'$. Choose $p\ge1$ such that $\varphi^n(K')\cap K'=\emptyset$ for all $n\ge p$.

Since $\varphi$ is continuous, injective and open, each $\varphi^{n}$ is a homeomorphism of $X$ onto the open set $O_n=\varphi^{n}(X)$; put $O_0=X$. For $j\in\Z$ define $h_j\colon X\to E$ and $S_j\subseteq X$ by
\[
h_j=\begin{cases} h\circ(\varphi^{pj})^{-1}\ \text{on } O_{pj},\quad 0 \text{ on } X\setminus O_{pj}, & j\ge0,\\[2pt] h\circ\varphi^{p|j|}, & j<0,\end{cases}
\qquad
S_j=\begin{cases}\varphi^{pj}(K'), & j\ge0,\\ (\varphi^{p|j|})^{-1}(K'), & j<0.\end{cases}
\]
Each $h_j$ vanishes outside $S_j$. Each $h_j$ is continuous: for $j<0$ this is clear; for $j\ge1$, $h\circ(\varphi^{pj})^{-1}$ is continuous on the open set $O_{pj}$ and vanishes on $O_{pj}\setminus S_j$, where $S_j$ is compact, so $h_j$ is continuous on each of the open sets $O_{pj}$ and $X\setminus S_j$, which cover $X$.

\emph{The sets $S_j$ are pairwise disjoint.} Let $i<j$ and $y\in S_i\cap S_j$. If $0\le i$, then $y=\varphi^{pi}(a)=\varphi^{pj}(b)$ with $a,b\in K'$, and injectivity gives $a=\varphi^{p(j-i)}(b)\in K'\cap\varphi^{p(j-i)}(K')=\emptyset$. If $i<0\le j$, then $\varphi^{p|i|}(y)\in K'$ and $y=\varphi^{pj}(b)$ with $b\in K'$, so $\varphi^{p(|i|+j)}(b)\in K'\cap\varphi^{p(|i|+j)}(K')=\emptyset$. If $i<j<0$, then $z=\varphi^{p|j|}(y)\in K'$ and $\varphi^{p(|i|-|j|)}(z)=\varphi^{p|i|}(y)\in K'$, again impossible. In each case we used $\varphi^{n}(K')\cap K'=\emptyset$ for $n\ge p$.

\emph{The family $(S_j)$ is locally finite.} Let $x\in X$ and let $L$ be a compact neighbourhood of $x$. Put $C=K'\cup L$ and choose $N$ with $\varphi^n(C)\cap C=\emptyset$ for $n\ge N$. If $j\ge1$ and $S_j\cap L\neq\emptyset$, then $\varphi^{pj}(C)\cap C\neq\emptyset$, so $pj<N$; if $j<0$ and $S_j\cap L\neq\emptyset$, then $\varphi^{p|j|}(L)\cap K'\neq\emptyset$, so $p|j|<N$. Hence only finitely many $S_j$ meet $L$.

Consequently $F=\sum_{j\in\Z}h_j$ is, near each point, a finite sum of continuous functions, so $F\in C(X,E)$. For $x\in K\subseteq K'=S_0$ we have $x\notin S_j$ for $j\neq0$, so $F(x)=h(x)=g(x)$.

\emph{Periodicity.} We claim $h_j(\varphi^p(y))=h_{j-1}(y)$ for all $j\in\Z$, $y\in X$. For $j\le0$, $h_j(\varphi^p y)=h(\varphi^{p(|j|+1)}y)=h_{j-1}(y)$. For $j\ge1$: if $y\in O_{p(j-1)}$, say $y=\varphi^{p(j-1)}(z)$, then $\varphi^p(y)=\varphi^{pj}(z)\in O_{pj}$ and $h_j(\varphi^py)=h(z)=h_{j-1}(y)$; if $y\notin O_{p(j-1)}$ (possible only when $j\ge2$), then $\varphi^p(y)\notin O_{pj}$ by injectivity of $\varphi^p$, and both sides vanish. Summing over $j$ gives $F\circ\varphi^p=F$.
\end{proof}

\begin{lemma}\label{lem:interval}
Let $J\subseteq\R$ be an open interval and $\psi\colon J\to J$ continuous without fixed points. Then $\psi$ is strongly run-away.
\end{lemma}

\begin{proof}
By the intermediate value theorem $\psi(x)-x$ has constant sign; say $\psi(x)>x$ for all $x$ (the other case is symmetric). Then every orbit $(\psi^n(x))_n$ is increasing. Let $K\subseteq[c,d]\subseteq J$ and $\delta=\min_{[c,d]}(\psi(x)-x)>0$, and let $N>(d-c)/\delta$. For $x\in[c,d]$, as long as $\psi^i(x)\in[c,d]$ for $i<n$, we have $\psi^n(x)\ge x+n\delta$; hence $\psi^i(x)>d$ for some $i\le N$, and then $\psi^n(x)>d$ for all $n\ge i$. Thus $\psi^n(K)\cap K=\emptyset$ for all $n\ge N$.
\end{proof}

\begin{lemma}\label{lem:semiconj}
Let $\theta\colon X\to Y$ be continuous, $\varphi\colon X\to X$, $\psi\colon Y\to Y$ with $\theta\circ\varphi=\psi\circ\theta$. If $\psi$ is (strongly) run-away, then so is $\varphi$.
\end{lemma}

\begin{proof}
If $K\subseteq X$ is compact, so is $\theta(K)$, and $\theta(\varphi^n(K))=\psi^n(\theta(K))$. If $x\in\varphi^n(K)\cap K$ then $\theta(x)\in\psi^n(\theta(K))\cap\theta(K)$.
\end{proof}

\begin{proof}[Proof of Corollary~\ref{cor:manifold}]
(b)$\Rightarrow$(a) is trivial and (a)$\Rightarrow$(c) is Theorem~\ref{thm:single}(3). For (c)$\Rightarrow$(b) it suffices, by Theorem~\ref{thm:chaos}, to show that $\varphi$ is open. Let $x\in X$ and choose charts $U\ni x$, $U'\ni\varphi(x)$ homeomorphic to $\R^d$. Then $\varphi$ restricts to a continuous injective map from the open set $U\cap\varphi^{-1}(U')\subseteq U\cong\R^d$ into $U'\cong\R^d$, which is open by Brouwer's invariance of domain~\cite[Thm.~2B.3]{Hatcher}. A second countable manifold is a separable locally compact metrizable space, so Theorem~\ref{thm:chaos} applies.
\end{proof}

\begin{proof}[Proof of Corollary~\ref{cor:interval}]
(b)$\Rightarrow$(a) and (c)$\Rightarrow$(a) are trivial. (a)$\Rightarrow$(d): by Theorem~\ref{thm:single}(1), $\varphi$ is injective and run-away, and a run-away map has no fixed points. (d)$\Rightarrow$(b),(c): $\varphi$ is strongly run-away by Lemma~\ref{lem:interval}, and open since an injective continuous map of an interval is strictly monotone; apply Theorem~\ref{thm:chaos}. The statement about hypercyclicity is Theorem~\ref{thm:single}(2).
\end{proof}

\begin{proof}[Proof of Corollary~\ref{cor:semiconj}]
By Lemmas~\ref{lem:interval} and~\ref{lem:semiconj}, $\varphi$ is strongly run-away; apply Theorem~\ref{thm:chaos}.
\end{proof}

\begin{example}
(a) Let $Z$ be a separable locally compact metric space, $J$ an open interval, $\psi\colon J\to J$ continuous, injective and without fixed points, and $\chi\colon Z\to Z$ an arbitrary homeomorphism. Then $\varphi=\psi\times\chi$ on $X=J\times Z$ is continuous, injective and open, and $\theta(t,z)=t$ semiconjugates $\varphi$ to $\psi$. By Corollary~\ref{cor:semiconj}, $C_\varphi$ is mixing and Devaney chaotic on $C(J\times Z,E)$. For $\chi=\mathrm{id}$ this is the example given in the original version of this paper.

(b) On $\R^d$, every translation $x\mapsto x+v$, $v\neq0$, and more generally every injective continuous $\varphi$ with $\langle\varphi(x)-x,v\rangle\ge c>0$ for all $x$ and some fixed $v$, is strongly run-away: with $\theta(x)=\langle x,v\rangle$ we have $\theta(\varphi^n x)\ge\theta(x)+nc$, so if $\theta(K)\subseteq[\alpha,\beta]$ then $\varphi^n(K)\cap K=\emptyset$ for $n>(\beta-\alpha)/c$. By Corollary~\ref{cor:manifold}, $C_\varphi$ is mixing and chaotic on $C(\R^d,E)$.
\end{example}

\begin{remark}[Openness]
Openness is not necessary for chaos. On $X=[0,\infty)$ the map $\varphi(x)=x+1$ is injective and strongly run-away but not open. Nevertheless $C_\varphi$ is chaotic on $C([0,\infty),E)$: given $g$ and $a>0$, choose an integer $p>a$, let $F=g$ on $[0,a]$, interpolate linearly on $[a,p]$ from $g(a)$ to $g(0)$, and extend $F$ $p$-periodically; then $F$ is continuous, $F\circ\varphi^p=F$ and $F=g$ on $[0,a]$.
\end{remark}

\section{Further remarks and questions}\label{sec:remarks}

\begin{remark}[Sensitivity]
Since $C(X,E)$ is an infinite metric space (with the metric $\rho$ of Section~\ref{sec:prelim}), the theorem of Banks et al.~\cite{Banks} shows that every Devaney chaotic $C_\varphi$ in Theorem~\ref{thm:chaos} also has sensitive dependence on initial conditions with respect to $\rho$. For semigroup actions, sensitivity does not in general follow from transitivity and density of periodic points, and sensitivity is a metric rather than a topological notion; for this reason we do not discuss Devaney chaos for semigroups of composition operators.
\end{remark}

\begin{remark}[More general targets]
The proofs use only that $E$ is a metrizable locally convex space (Dugundji's theorem, and a countable dense set for hypercyclicity). Replacing $\norm{\cdot}$ by an increasing sequence of seminorms, Theorems~\ref{thm:main}, \ref{thm:single} and~\ref{thm:chaos} extend to separable Fr\'echet spaces $E$.
\end{remark}

\begin{question}
Does condition \TE{} depend on $E$? That is, if $G$ satisfies (T$_\R$), does it satisfy (T$_E$) for every normed space $E$?
\end{question}

\begin{question}
Find a characterization of Devaney chaos for $C_\varphi$ on $C(X,E)$. In particular, is every injective strongly run-away map $\varphi$ on a separable locally compact metric space such that $C_\varphi$ has a dense set of periodic points?
\end{question}

\begin{question}
Give natural (e.g.\ manifold) examples of injective run-away maps that are not strongly run-away, i.e.\ of hypercyclic but non-mixing operators $C_\varphi$ on $C(X,E)$.
\end{question}

\section*{Acknowledgement}
I would like to thank the referee for making many useful suggestions and for pointing out the references~\cite{BBCP,CP,dLE,DSW}. This revised version is prepared by an AI assistant (Claude).

\end{document}